\documentclass[12pt]{amsart}
\usepackage[utf8]{inputenc}
\usepackage[margin=2.5cm]{geometry}
\usepackage{amssymb,mathtools,mathrsfs,bbm}
\usepackage{enumerate}
\usepackage{graphicx,titletoc,xcolor}
\definecolor{darkblue}{rgb}{0.0, 0.0, 0.55}
\definecolor{bordeaux}{rgb}{0.34, 0.01, 0.1}
\usepackage[colorlinks,linkcolor=bordeaux,citecolor=darkblue,urlcolor=black,hypertexnames=true]{hyperref}
\usepackage{tikz-cd}
\usepackage{euflag}

\makeatletter
\def\moverlay{\mathpalette\mov@rlay}
\def\mov@rlay#1#2{\leavevmode\vtop{
		\baselineskip\z@skip \lineskiplimit-\maxdimen
		\ialign{\hfil$#1##$\hfil\cr#2\crcr}}}
\makeatother
\newcommand{\lfre}{{<}}
\newcommand{\rfre}{{>}}
\newcommand{\lfree}{\moverlay{(\cr<}}
\newcommand{\rfree}{\moverlay{)\cr>}}
\newcommand{\kk}{{\mathbbm k}}
\newcommand{\mx}[1]{\lfre #1\rfre}
\newcommand{\px}[1]{\kk\lfre #1\rfre}
\newcommand{\sx}[1]{\kk\lfre\!\!\!\lfre #1\rfre\!\!\!\rfre}
\newcommand{\rx}[1]{\kk\lfree #1\rfree}

\DeclareMathOperator{\GL}{GL}

\DeclareMathOperator{\dom}{dom}

\DeclareMathOperator{\im}{im}
\DeclareMathOperator{\tr}{tr}

\newcommand{\N}{{\mathbb N}}
\newcommand{\Z}{{\mathbb Z}}
\newcommand{\C}{{\mathbb C}}

\newcommand{\cO}{{\mathcal O}}

\newcommand{\gen}{{\rm g}}

\newtheorem{theorem}{Theorem}[section]
\newtheorem{proposition}[theorem]{Proposition}
\newtheorem{lemma}[theorem]{Lemma}
\newtheorem{corollary}[theorem]{Corollary}

\makeatletter
\newtheorem*{rep@thm}{\rep@title}
\newcommand{\newreptheorem}[2]{%
	\newenvironment{rep#1}[1]{%
		\def\rep@title{#2 \ref{##1}}%
		\begin{rep@thm}}%
		{\end{rep@thm}}}
\makeatother
\newreptheorem{thm}{Theorem}
\newreptheorem{cor}{Corollary}

\newtheorem{mainthm}{Theorem}

\newtheorem{maincor}[mainthm]{Corollary}

\theoremstyle{definition}

\newtheorem{remark}[theorem]{Remark}

\title{Commutators in free skew fields}

\author[Igor Klep]{Igor Klep\textsuperscript{1}}
\address{Faculty of Mathematics and Physics, Department of Mathematics, University of Ljubljana \& FAMNIT, University of Primorska \& Institute of Mathematics, Physics and Mechanics, Ljubljana, Slovenia}
\email{igor.klep@fmf.uni-lj.si}
\thanks{\textsuperscript{1} Supported by the Slovenian Research Agency 
	program P1-0222 and grants J1-50002,
	N1-0217,
J1-60011, J1-50001, J1-3004 and J1-60025. Part of this work was performed within the project COMPUTE, funded within the QuantERA II 
Programme that has received funding from the EU's H2020 research and innovation programme under the GA No.~101017733 {\normalsize\euflag}. Part of this work was performed within the project SDPcode, funded within the QuantERA III 
Programme that has received funding from the EU's Horizon research and innovation programme under the GA No.~101212998 {\normalsize\euflag}.}

\author[Jurij Vol\v{c}i\v{c}]{Jurij Vol\v{c}i\v{c}\textsuperscript{2}}
\address{Department of Mathematics, University of Auckland, New Zealand}
\email{jurij.volcic@auckland.ac.nz}
\thanks{\textsuperscript{2} Supported by the Marsden Fund MFP-UOA2528 from Government funding, administered by the Royal Society Te Ap\=arangi, and the New Staff grant 3734441 by the Faculty of Science, the University of Auckland. 
Supported in part by a grant of access to OpenAI models through the ChatGPT for Academic Researchers program.}

\date{\today}



\begin{document}

\begin{abstract}
We characterise additive and multiplicative commutators in free
skew fields in terms of traces and determinants of matrix
evaluations. Namely, a noncommutative rational function is a sum
of additive commutators if and only if all its matrix evaluations
have trace zero, and it is a product of multiplicative commutators
if and only if all its matrix evaluations have determinant one.
As applications, we characterise pairs of noncommutative rational
functions with equal operator or Frobenius norms, or with
identical spectra. 
Finally, using valuations and free
Poisson fields, we show that the identity and the generators of a free skew field
cannot be expressed as a nonempty product of additive
commutators, answering a question of Gardella and Thiel
in the negative.
\end{abstract}

\maketitle

\tableofcontents

\section{Introduction}\label{s:intro}

Commutators measure the failure of commutativity, and their relationship
with matrix invariants has been studied since the early days of matrix
algebra. Shoda \cite{Sho37} proved in characteristic 0 that a matrix
is an additive commutator if and only if its trace vanishes;
Albert and Muckenhoupt \cite{AM57} extended this to arbitrary fields; see also \cite{Suw23} for more recent
work over commutative rings. 
For finite-dimensional simple algebras, Amitsur and Rowen
\cite{AR94} proved that every element of reduced trace zero
is a sum of at most two additive commutators. The distinction
between a single commutator and a sum of commutators is
illustrated by the examples of Rosset and Rosset \cite{RR00}
and, more recently, Jang \cite{Jang26}.
In the free algebra, 
a noncommutative polynomial
has zero trace at matrices of every size precisely when it is a sum of
commutators. This observation appears in 
\cite{KS08,BK}, and underlies the tracial
Nullstellensatz of \cite{KSp14}. It is also closely
related to the study of values of noncommutative polynomials and trace
positivity; see, for instance, \cite{Bre20,KSV23}.
The invariant theory and trace identities of generic matrices
were developed by Procesi \cite{Pro76} and Razmyslov
\cite{Raz74}, while rational identities in division algebras
were investigated by Bergman \cite{Ber76}.

On the multiplicative side, the Dieudonn\'e determinant
\cite{Die43} and the theory of universal fields of fractions developed
by Cohn \cite{Coh} provide a natural framework for commutator subgroups
of division rings. Amitsur's theory of rational identities \cite{Ami}
and the matricial approach to noncommutative rational functions of
Helton, McCullough and Vinnikov \cite{HMV06}, and
Kalyuzhnyi-Verbovetskyi and Vinnikov \cite{KVV1,KVV2} allow one to
translate questions in a free skew field into questions about matrix
evaluations. The connection between factorisation in the free algebra
and irreducibility of generic determinants was developed in \cite{HKV18,HKV}; generic spectral properties
of rational evaluations were subsequently studied in \cite{BV}. Questions about pointwise rank,
spectral and norm equivalences of noncommutative polynomials are
considered in \cite{SV}.

A different problem concerns \emph{products of additive commutators}.
Building on earlier work of Baxter \cite{Bax65}, Gardella and Thiel
\cite{GT} showed that every element of a noncommutative division ring
is a sum of two products of pairs of additive commutators, and asked
whether every element must itself be a product of two such commutators
\cite[Question~5.9]{GT}. The answer is positive for skew Laurent series
division rings by Jang and Ke \cite{JK26}, and for finite-dimensional
noncommutative division algebras by Bre\v{s}ar, Jang and Robert
\cite{BJR}. Furthermore, Bre\v{s}ar, Gardella and Thiel \cite{BGT}
proved that every element of $D^{n\times n}$, for $n\ge2$ and a division ring
$D$ with infinite centre, is a product of two additive commutators.
In contrast, we show that even the identity element of a free skew
field cannot be expressed as a nonempty product of additive
commutators. The proof exploits a valuation due to Lichtman \cite{Lic}
and algebraic dependence in free Poisson fields as studied by
Makar-Limanov and Shestakov \cite{MS}.

\subsection{Notation}\label{ss:intro-notation}
Throughout, $\kk$ is a field of characteristic 0 and
$x=(x_1,\ldots,x_\gen)$ is a tuple of freely noncommuting variables.
We write
\[
 R=\px{x_1,\ldots,x_\gen},\qquad
 U=\rx{x_1,\ldots,x_\gen}
\]
for the free algebra and its universal skew field of fractions,
respectively. The construction and universal property of $U$ go back
to Cohn; see \cite[Chapter~7]{Coh}. Its elements can be viewed as
\emph{noncommutative rational functions}
\cite[Section~2]{KVV1}; see also \cite{Ami,HMV06,KVV2,VolReal18} for
rational identities, the free-function viewpoint and realisation
theory. 
Noncommutative rational series and their realisations also arise
in multidimensional systems theory \cite{BGM}, while
linearisations of noncommutative rational functions play an
important role in free probability and random matrix theory
\cite{HMS,CMMPY}.
The domain of a noncommutative rational function $r$ is denoted $\dom r$ and consists of all $\gen$-tuples of square matrices over $\kk$ at which $r$ can be evaluated in a natural way. 
Precise definitions and conventions used in the proofs are
recalled in Section~\ref{s:prelim}.

For a $\kk$-algebra
$A$, let $[A,A]$ be the $\kk$-linear span of the additive commutators
$[a,b]=ab-ba$ for $a,b\in A$. For a group $G$, let $G'$ denote its commutator
subgroup. In particular, $(U^\times)'$ consists of finite products
of multiplicative commutators $aba^{-1}b^{-1}$ for
$a,b\in U^\times$. When $\kk=\C$, the operator and Frobenius norms
are denoted by $\|\cdot\|$ and $\|\cdot\|_F$, respectively.

\subsection{Main results}\label{ss:intro-results}

Our first result characterises additive and multiplicative
commutators in the free skew field by traces and determinants of
matrix evaluations.

\begin{mainthm}\label{thm:mainA}
Let $r\in U$. The following hold:
\begin{enumerate}[\rm (a)]
\item $\tr r(X)=0$ for every $X\in\dom r$ if and only if
$r\in[U,U]$.
\item $\det r(X)=1$ for every $X\in\dom r$ if and only if
$r\in(U^\times)'$.
\end{enumerate}
\end{mainthm}

Parts (a) and (b) are proved as Theorems~\ref{t:add} and
\ref{t:mult}, respectively. The proof of the former in
Section~\ref{s:add} is based on a commutator criterion for rational
formal power series (Proposition~\ref{p:series}), which in turn relies on a lifting
argument from cyclic homology (Lemma~\ref{l:ch}). Part (b) is
proved in Section~\ref{s:mult} using Cohn's description of the abelianisation of $U^\times$ (Proposition~\ref{p:coh}) and
factorisation of generic determinants \cite{Coh,HKV}.

Theorem \ref{thm:mainA} addresses the structure behind two eigenvalue relations for noncommutative rational functions that hold for evaluations on matrices of all sufficiently large sizes. 
It is an open question whether a nonconstant noncommutative rational function can admit two or more algebraically independent eigenvalue relations that persist through all large sizes (cf. the end of \cite[Section 5.2]{DV}). A noncommutative rational function can also satisfy an eigenvalue relation that is not an obvious consequence of trace or determinant constraints (Remark \ref{r:semaev}).

As consequences of Theorem \ref{thm:mainA}, we obtain rigidity statements for norms and
spectra of noncommutative rational functions, 
extending the polynomial results of \cite{SV}.

\begin{maincor}\label{cor:mainB}
Let $\kk=\C$ and $r,s\in U$.
\begin{enumerate}[\rm (a)]
\item $\|r(X)\|_F=\|s(X)\|_F$ for every
$X\in\dom r\cap \dom s$ if and only if
\[
r=\alpha+u,\qquad s=\xi\alpha+\zeta u
\]
for some $\alpha\in\C$, $u\in U$ and $\xi,\zeta\in\C$ of
modulus one such that either $\xi=\zeta$ or $u\in[U,U]$.
\item $\|r(X)\|=\|s(X)\|$  for every $X\in\dom r\cap\dom s$ if and only if $s=\zeta r$ for some $\zeta\in\C$ of modulus one.
\item The spectra of $r(X)$ and $s(X)$ coincide for every
$X\in\dom r\cap \dom s$ if and only if $s=uru^{-1}$ for some
$u\in U^\times$.
\end{enumerate}
\end{maincor}

Parts (a) and (b) are Corollary~\ref{c:norm}, proved in
Section~\ref{s:add}, and part (c) is Corollary~\ref{c:isospec}, proved in
Section~\ref{s:mult}. 

Finally, we answer \cite[Question~5.9]{GT} in the negative,
in contrast with positive results for other classes of division
rings \cite{JK26,BJR}.

\begin{mainthm}\label{thm:mainC}
$1$ and the $x_i$ are not nonempty products of additive commutators in $U$.
\end{mainthm}

Theorem~\ref{thm:mainC} is proved as
Corollary~\ref{c:one} in Section~\ref{s:one}, using
Lichtman's valuation \cite{Lic} and an obstruction in a free
Poisson field \cite{MS}. It contrasts also with the matrix-ring
result of Bre\v{s}ar, Gardella and Thiel
\cite[Theorem~4.4]{BGT}, which applies to $U^{n\times n}$ for $n\ge2$.

\subsection{Use of artificial intelligence}\label{ss:intro-ai}
OpenAI's 
GPT-5.6 Sol and GPT-6 Astra were
used during the preparation of this manuscript 
to draft initial proof arguments following the authors' guidance and strategy, 
locate existing results, 
assist with the literature review, and check the soundness of the authors' final proofs. 
The authors take responsibility for the
correctness and originality of the results, the proofs, and the
bibliographic information.

\section{Preliminaries}\label{s:prelim}

Throughout the paper let $\kk$ be a field of characteristic 0. 
All the rings and ring homomorphisms are unital. 
Let $R=\px{x_1,\dots,x_\gen}$ be the free $\kk$-algebra generated by $x_1,\dots,x_\gen$, and let $U=\rx{x_1,\dots,x_\gen}$ be its universal skew field of fractions \cite[Chapter 7]{Coh}, the free skew $\kk$-field generated by $x_1,\dots,x_\gen$. 
An $m\times m$ matrix over $R$ is \emph{full} if it does not factor over $R$ as $AB$ where $A$ has fewer than $m$ columns. Then, $U$ is characterised as follows: $R\subseteq U$, $U$ is generated by $R$ as a skew field, and every full matrix over $R$ is invertible over $U$ \cite[Section 7.5]{Coh}.

By \cite[Section 2]{KVV1} (see also \cite{Ami}),  we can view elements of $U$ as equivalence classes of formal rational expressions in $x_1,\dots,x_\gen$ over $\kk$. More precisely, given a formal rational expression $\mathfrak r$, let $\dom_n\mathfrak r\subseteq (\kk^{n\times n})^\gen$ be the set of $n\times n$ matrix tuples on which $\mathfrak r$ can be evaluated. 
Note that $\dom_n\mathfrak r$ is a Zariski open subset of $(\kk^{n\times n})^\gen$. We say that $\mathfrak r$ is \emph{nondegenerate} if $\dom\mathfrak r=\bigsqcup_{n\in\N}\dom_n\mathfrak r$ is nonempty. Then, two nondegenerate expressions $\mathfrak r_1$ and $\mathfrak r_2$ are equivalent if $\mathfrak r_1(X)=\mathfrak r_2(X)$ for all $X\in\dom \mathfrak r_1\cap\dom\mathfrak r_2$. Equivalence classes of nondegenerate expressions are \emph{noncommutative rational functions}, and they form a skew field isomorphic to $U$ by \cite[Proposition 2.2]{KVV1}. 
We define the \emph{domain} of $r\in U$ as
$$\dom r= \bigcup_{\mathfrak r\in r}\dom\mathfrak r,\qquad
\dom_n r=\dom r\cap \bigl(\kk^{n\times n}\bigr)^\gen.$$
The set $\dom _n r$ is Zariski open in the affine space $(\kk^{n\times n})^\gen$, and invariant under simultaneous conjugation by invertible matrices. Furthermore, $(\dom_n r)^{\oplus m}\subseteq \dom_{mn} r$ for all $m,n\in\N$. 
Evaluations of $r$ respect simultaneous conjugations and direct sums wherever they are defined \cite{KVV1,KVV2}. 
Consequently, many properties of a noncommutative rational function only need to be checked on $\dom_n r$ for infinitely many $n$. For example, if $r$ and $s$ have the same trace on $\dom_n r\cap\dom_n s$ for infinitely many $n$, then they have the same trace on $\dom r\cap \dom s$. 
We use such conclusions tacitly throughout the paper. 
For a detailed treatment of the distinction between domains of rational expressions and rational functions, and their descriptions by minimal realisations, see 
\cite{VolDom17,VolReal18,PV}.

Given a field extension $K$ of $\kk$, denote 
$U_K=K\lfree x_1,\dots,x_\gen \rfree$. 
If $K$ is a finite extension of $\kk$, then $U_K\cong K\otimes_\kk U$ \cite[Exercise 7.5.13]{Coh}. 
If $t$ is transcendental over $\kk$, then $U_{\kk(t)}\cong U(t)$, the classical ring of quotients of the Ore domain $U[t]$ \cite[Exercise 7.5.14]{Coh}.

\section{Additive commutators}\label{s:add}

In this section, we characterise additive commutators in the free
skew field in terms of traces of matrix evaluations
(Theorem~\ref{t:add}). As an application, we describe pairs
of noncommutative rational functions whose matrix evaluations
have equal operator or Frobenius norms (Corollary~\ref{c:norm}).

\def\rat{{\operatorname{rat}}}
\def\cyc{\overset{{\rm\scriptscriptstyle cyc}}{\sim}}

We start with an auxiliary homological lemma (see \cite[Chapter 9]{Wei} for background on Hochschild and cyclic homology needed in the proof).

\begin{lemma}\label{l:ch}
Let $A$ be a finite-dimensional $\kk$-algebra, with a bilinear form $\beta:A\times A\to \kk$ satisfying
\[
 \beta(a,b)=-\beta(b,a),\quad
 \beta(ab,c)=\beta(a,bc)+\beta(b,ca)
 \qquad \text{for all }a,b,c\in A.
\]
There are a finite-dimensional $\kk$-algebra $E$, a surjective homomorphism $\pi:E\to A$, and a linear functional $\lambda:E\to \kk$ such that
\[
 \beta(\pi(a),\pi(b))=\lambda([a,b])
 \qquad \text{for all }a,b\in E.
\]
\end{lemma}

\begin{proof}
Let $J$ be the Jacobson radical of $A$. Then, $A/J$ is semisimple and $J$ is a nilpotent ideal. By Wedderburn's principal theorem \cite[Exercise 9.3.1]{Wei}, there is a vector space decomposition $A=S\oplus J$ where $S$ is a subalgebra of $A$ isomorphic to $A/J$. 
Consider the graded tensor algebra
$$T=\bigoplus_{n=0}^\infty T_n\coloneq 
\bigoplus_{n=0}^\infty J^{\otimes_S n}=S\oplus J\oplus (J\otimes_S J)\oplus (J\otimes_S J\otimes_S J)\oplus\cdots$$
and let $m:T\to A$ be the multiplication homomorphism. 

By \cite[Lemma~4.2 and Remark~4.3]{Cib} (see also \cite[Proposition 9.1.6]{Wei}), there is a projective $T$-bimodule resolution
\[
 0\longrightarrow T\otimes_S J\otimes_S T
 \longrightarrow T\otimes_S T\longrightarrow T\longrightarrow0 .
\]
Thus, the second Hochschild homology of $T$ vanishes. Since $S$ is a semisimple finite-dimensional algebra of characteristic 0, the first cyclic homology of $S$ vanishes by \cite[Example 9.6.14]{Wei}. 
Thus, the first cyclic homology of $T$ vanishes by Goodwillie's theorem \cite[Theorem 9.9.1]{Wei}.
Denote
\begin{align*}
C_1&=(T\otimes_\kk T)\big/ \operatorname{span}_\kk\{t_0\otimes t_1+t_1\otimes t_0\colon t_i\in T\}, \\
C_2&=(T\otimes_\kk T\otimes_\kk T)\big/ \operatorname{span}_\kk\{t_0\otimes t_1\otimes t_2-t_2\otimes t_0\otimes t_1\colon t_i\in T\}, \\
\partial_1&:C_1\longrightarrow T, \quad\ \partial_1(t_0\otimes t_1)=[t_0,t_1], \\
\partial_2&:C_2\longrightarrow C_1, \quad \partial_2(t_0\otimes t_1\otimes t_2)=t_0t_1\otimes t_2-t_0\otimes t_1t_2+t_2t_0\otimes t_1.
\end{align*}
Vanishing of the first cyclic homology means $\im \partial_2=\ker \partial_1$ by \cite[Lemma 9.6.10]{Wei}. Taking $\kk$-duals gives $\im \partial_1^*=\ker \partial_2^*$. 

Since the bilinear form $\beta$ is antisymmetric, there is $\widetilde\beta\in C_1^*$ determined by $\widetilde\beta(t_0\otimes t_1)=\beta(m(t_0),m(t_1))$. Furthermore, $\partial_2^*(\widetilde\beta)=0$ by the cocycle condition on $\beta$. Since $\ker \partial_2^*= \im \partial_1^*$, there is $\Lambda\in T^*$ such that $\widetilde\beta=\partial_1^*(\Lambda)$, and so
\begin{equation}\label{e:preL}
 \beta\big(m(a),m(b)\big)=\Lambda([a,b])\qquad 
 \text{for all }a,b\in T.
\end{equation}

Since $J$ is nilpotent, there is $p\in\N$ such that $J^p=0$. Since $m(T_n)=0$ for $n\geq p$, we have
$\widetilde\beta(T_m,T_n)=0$ whenever $m+n\ge 2p-1$. 
Consequently, we can modify the functional $\Lambda$ so that it vanishes on $\bigoplus_{n\ge2p-1}T_n$, and \eqref{e:preL} still holds. 
Then, both $\Lambda$ and $m$ factor through the finite-dimensional algebra
$$E= T\Big/\bigoplus_{n\geq 2p-1}T_n,$$
giving the desired $\lambda$ and $\pi$.
\end{proof}

We equip the free monoid $W=\mx{x_1,\dots,x_\gen}$ with the standard length $|w|\in\N_0$ for $w\in W$. 
Let $\widehat{R}=\sx{x_1,\dots,x_\gen}$ be the $\kk$-algebra of formal power series. For $s\in\widehat{R}$, let $s_w\in\kk$ denote its coefficient at $w\in W$; 
that is, $s=\sum_{w\in W}s_w w$.
The rational closure of $R$ in $\widehat{R}$ is the  subalgebra of rational series $\widehat{R}_\rat$. There is a natural isomorphism $\widehat{R}_\rat\cong \{r\in U\colon 0\in\dom r\}$. 
By Sch\"utzenberger's theorem \cite[Theorem 7.1]{BR}, every $r\in\widehat{R}_\rat$ admits a \emph{linear representation},
$$r=c^\intercal(I-x_1A_1-\cdots -x_\gen A_\gen)^{-1}b$$
where $b,c\in \kk^{n\times 1}$ and $A_i\in \kk^{n\times n}$. Conversely, such a linear representation determines a rational series in $\widehat{R}$.

We say that $s=\sum_{w\in W} s_w w\in \widehat{R}$ is \emph{cyclically equivalent to 0} \cite[Definition 4.1 and Remark 4.2]{BK} if $s_1=0$ and for every $w=x_{i_1}\cdots x_{i_\ell}$, $\ell\in\N$,
$$\sum_{k=1}^{\ell}s_{
x_{i_k}\cdots x_{i_\ell}x_{i_1}\cdots x_{i_{k-1}}
}=0.$$

\begin{lemma}\label{l:cyc}
Let $r\in\widehat{R}_\rat$. 
If $\tr r=0$ on $\dom r$, then $r$ is cyclically equivalent to 0.
\end{lemma}

\begin{proof}
Let $X\in(\kk^{n\times n})^\gen$ be arbitrary. Since $0\in\dom r$, it follows that $r$ is defined at $tX\in (\kk[\![t]\!]^{n\times n})^\gen$, and so
$$r(tX)=\sum_w r_w t^{|w|} w(X).$$
Then, $\tr r(tX)=0$ implies 
$\sum_{|w|=\ell} r_w \tr w(X)=0$ 
for all $\ell\in\N_0$. Since $X$ was arbitrary, the homogeneous polynomial 
$\sum_{|w|=\ell} r_w w\in R$ is cyclically equivalent to 0 by \cite[Corollary 4.8]{BK}, for every $\ell\in\N_0$. Consequently, $r$ is cyclically equivalent to 0.
\end{proof}

Next, we show that if a rational series is cyclically equivalent to 0, then it is a sum of commutators of rational series.

\begin{proposition}\label{p:series}
Let $r\in\widehat{R}_\rat$ be cyclically equivalent to $0$. There are $r_1,\dots,r_\gen\in \widehat{R}_\rat$ such that
$$r=[x_1,r_1]+\cdots+[x_\gen, r_\gen].$$
\end{proposition}

\begin{proof}
Choose a linear representation $r=c^\intercal (I-\sum_ix_iA_i)^{-1}b$. Let $\mu:R\to \kk^{n\times n}$ be the unital $\kk$-algebra homomorphism determined by $\mu(x_i)=A_i$, and let $C=bc^\intercal \in\kk^{n\times n}$. Then, $r_w=\tr(C\mu(w))$ for all $w\in W$. 

Let $\delta:R\to \kk^{n\times n}$ be the $\mu$-derivation determined by 
$\delta(x_i)=CA_i$ for $i=1,\dots,\gen$; 
that is, $\delta(fg)=\mu(f)\delta(g)+\delta(f)\mu(g)$ for all $f,g\in R$.
For every $w=x_{i_1}\cdots x_{i_\ell}$,
$$\tr \delta(w)=
\sum_{k=1}^{\ell} \tr\bigl(CA_{i_k}\cdots A_{i_\ell}A_{i_1}\cdots A_{i_{k-1}}\bigr)
=\sum_{k=1}^{\ell}r_{
x_{i_k}\cdots x_{i_\ell}x_{i_1}\cdots x_{i_{k-1}}
}=0$$
by the assumption on $r$.
Therefore,
\begin{equation}\label{e:tr0}
 \tr\delta(f)=0 \qquad \text{for all }f\in R.
\end{equation}
Let $A=\mu(R)$. The rule
$$\beta(\mu(f),\mu(g))\coloneq \tr\bigl(\delta(f)\mu(g)\bigr)$$
defines a bilinear form on $A$. Indeed, if $\mu(f)=0$, then 
\[
 \tr\bigl(\delta(f)\mu(g)\bigr)=\tr\delta(fg)=0
\]
by the Leibniz rule for $\delta$ and \eqref{e:tr0}. The latter two also imply
\[
 \beta(\mu(f),\mu(g))+\beta(\mu(g),\mu(f))
   =\tr\delta(fg)=0,
\]
so $\beta(a,b)=-\beta(b,a)$.
Furthermore, the Leibniz rule and cyclicity of trace give
$\beta(ab,c)=\beta(a,bc)+\beta(b,ca)$.
By Lemma~\ref{l:ch}, there is an algebra $E\subseteq \kk^{N\times N}$ with a surjective homomorphism $\pi:E\to A$ and a functional $\lambda:E\to \kk$ such that $\beta(\pi(a),\pi(b))=\lambda([a,b])$ for all $a,b\in E$.

Choose a unital homomorphism $\nu:R\to E$ such that $\mu=\pi\circ\nu$, and define
$$r_i\coloneq-\sum_{w\in W}\lambda\bigl(\nu(wx_i)\bigr) w\in \widehat{R}$$
for $i=1,\dots,\gen$.
Note that $r_i$ is a rational series. Indeed, $X\mapsto \lambda(X \nu(x_i))$ is a functional on $E\subseteq\kk^{N\times N}$, so there are $b_{i1},c_{i1},\dots, b_{iN},c_{iN}\in \kk^{N\times 1}$ such that this functional agrees with $X\mapsto \sum_j c_{ij}^\intercal X b_{ij}$ on $E$; then,
$$r_i=-\sum_{j=1}^N c_{ij}^\intercal 
\big(I-x_1\nu(x_1)-\cdots-x_\gen \nu(x_\gen)\big)^{-1} b_{ij}\in \widehat{R}_\rat.$$

Let $w\in W$ be arbitrary. If $|w|\le 1$, then $r_w=0$ because $r$ is cyclically equivalent to 0; likewise, the coefficient of $w$ in $\sum_k[x_k,r_k]$ is 0. 
If $|w|\ge 2$, then $w=x_iux_j$, and the coefficient of $w$ in $\sum_k[x_k,r_k]$ equals
\begin{align*}
(r_i)_{ux_j}-(r_j)_{x_iu}
&= \lambda\bigl([\nu(x_i),\nu(ux_j)]\bigr) \\
&=\beta\bigl(\mu(x_i),\mu(ux_j)\bigr) \\
&=\tr\bigl(C\mu(x_i)\mu(ux_j)\bigr) \\
&= c^\intercal \mu(x_iux_j) b,
\end{align*}
which is the coefficient of $w$ in $r$. Thus, $r=\sum_k[x_k,r_k]$.
\end{proof}

Using expansions of noncommutative rational functions about matrix points and Morita invariance, Proposition \ref{p:series} leads to the following characterisation of noncommutative rational functions with trace zero.

\begin{theorem}\label{t:add}
Let $r\in U$. Then, $\tr r= 0$ on $\dom r$ if and only if $r\in[U,U]$.
\end{theorem}

\begin{proof}
$(\impliedby)$ If $r=[a_1,b_1]+\cdots +[a_\ell, b_\ell]$ for $a_i,b_i\in U$, then $\tr r= 0$ on $\bigcap_i\dom a_i\cap \dom b_i$, and so $\tr r= 0$ on $\dom r$ by 
eventual nonemptiness and Zariski density.

$(\implies)$ Assume $\tr r=0$ on $\dom r$, and fix $X\in\dom_m r$. 
Let $y_{k ij}$ for $1\le k\le \gen$ and $1\le i,j\le m$ be auxiliary free variables, and denote 
$$S=\sx{y_{k ij}\colon i,j,k}_{\rat},\qquad 
Y=\big( (y_{1 ij})_{i,j},\dots,(y_{\gen ij})_{i,j}\big)\in \big(S^{m\times m}\big)^\gen.$$
Then, the $m\times m$ rational matrix $M=r(X+Y)$ is defined at the origin, so we can view it as a matrix over $S$. 
Since $\tr r=0$ on $\dom r$, it follows that $\tr M=0$ on $\dom M$. 
Denote $T=M_{11}+\cdots+M_{mm}\in S$; then $\tr T=0$ on $\dom T$. Thus, $T\in [S,S]$ by Lemma \ref{l:cyc} and Proposition \ref{p:series}. The ordinary matrix trace induces the Morita invariance isomorphism \cite[Corollary 9.5.8]{Wei}
$$S^{m\times m}\big/ \big[S^{m\times m},S^{m\times m}\big]\cong S/ [S,S].$$
Hence, $M\in [S^{m\times m},S^{m\times m}]$.

Let $t$ be a central indeterminate. There is a homomorphism
$$S\longrightarrow U(t),\qquad 
y_{k ij}\longmapsto (1-t)
\bigl(\delta_{ij}x_k-(X_k)_{ij}\bigr).$$
Indeed, this rule is well-defined at $t=1$, and consequently at the generic $t$.
Applying this homomorphism entrywise to $M\in [S^{m\times m},S^{m\times m}]$ results in
\begin{equation}\label{e:Xt}
r\bigl(
tX+{(1-t)}(x_1I_m,\dots,x_\gen I_m)
\bigr)\in [U(t)^{m\times m},U(t)^{m\times m}].
\end{equation}
Let us embed $U(t)$ into $U(\!(t)\!)$; note that the constant term of a commutator in $U(\!(t)\!)^{m\times m}$ is a sum of commutators in $U^{m\times m}$. Thus, 
\eqref{e:Xt} implies $r I_m\in [U^{m\times m},U^{m\times m}]$. Applying the ordinary matrix trace gives $m r\in [U,U]$, so $r\in [U,U]$.
\end{proof}

The interplay between trace and commutators allows us to describe pairs of noncommutative rational functions whose values have the same norm. This is recorded in Corollary \ref{c:norm} below, which extends the polynomial statement \cite[Theorem 6.3]{SV}.

Let $\kk=\C$ be the field of complex numbers. It is well-known that 
\begin{equation}\label{e:dense}
\{(X,X^*)\colon X\in (\kk^{n\times n})^\gen\} 
\quad\text{is Zariski dense in}\quad
(\kk^{n\times n})^{2\gen} 
\quad\text{for every }n\in\N,
\end{equation}
where $X^*=(X_1^*,\dots,X_\gen^*)$. 
Let $x_1^*,\dots,x_\gen^*$ be the formal adjoints of $x_1,\dots,x_\gen$. The skew field $\rx{x_1,\dots,x_\gen,x_1^*,\dots,x_\gen^*}$ is endowed with a natural skew-linear involution that exchanges $x_j$ and $x_j^*$. By \eqref{e:dense}, every $r\in \rx{x_1,\dots,x_\gen,x_1^*,\dots,x_\gen^*}$ admits evaluations $r(X,X^*)$ for generic $X\in(\kk^{n\times n})^\gen$ and sufficiently large $n$, and all evaluations of this type determine $r$. Note that the free product of $\rx{x_1,\dots,x_\gen}$ and $\rx{x_1^*,\dots,x_\gen^*}$ over $\kk$ naturally embeds into $\rx{x_1,\dots,x_\gen,x_1^*,\dots,x_\gen^*}$, which has the following consequence.

\begin{lemma}\label{l:norm}
Let $\kk=\C$, $r\in U\setminus\kk$ and $s\in U^\times$.
\begin{enumerate}[\rm (a)]
\item If $rr^*-r^*r=ss^*-s^*s$, then $s=\zeta r+\lambda$ for $\lambda,\zeta\in\kk$ with $|\zeta|=1$.
\item If $(rr^*)^p-(r^*r)^p=(ss^*)^p-(s^*s)^p$ for $p=-1,1$, then $s=\zeta r$ for $\zeta\in\kk$ with $|\zeta|=1$.
\end{enumerate}
\end{lemma}

\begin{proof}
(a): Let $U=\kk\oplus U_0$ for some subspace $U_0$ of $U$, and write $r=\alpha+r_0$ and $s=\beta+s_0$ with respect to this decomposition. By the aforementioned natural embedding, we can consider the equation $rr^*-r^*r=ss^*-s^*s$ in
$$\kk\oplus U_0\oplus U_0^*\oplus 
(U_0\otimes U_0^*)\oplus (U_0^*\otimes U_0),$$
where $U_0^*\subset \rx{x_1^*,\dots,x_\gen^*}$ is the image of $U_0$ under the involution. Then, $r_0r_0^*=s_0s_0^*$ in $U_0\otimes U_0^*$, so $s_0=\zeta r_0$ for some $\zeta\in \kk$. Furthermore, $|\zeta|=1$ by $s_0s_0^*=r_0r_0^*$ since $r_0\neq0$.

(b): By (a) we have $s=\zeta r+\lambda$ and $s^{-1}=\zeta' r^{-1}+\lambda'$ for $\zeta,\lambda,\zeta',\lambda'\in\kk$ with $|\zeta|=|\zeta'|=1$. Then $r=\zeta' (\zeta r+\lambda)+\lambda' r(\zeta r+\lambda)$ is a quadratic equation in $r$, with constant coefficient $\zeta'\lambda$. Since a nonconstant element of a free skew field is transcendental over $\kk$, it follows that $\zeta'\lambda=0$, and so $\lambda=0$.
\end{proof}

Next, we can characterise pairs of noncommutative rational functions with the same Frobenius or operator norms.

\begin{corollary}\label{c:norm}
Let $\kk=\C$ and $r,s\in U$.
\begin{enumerate}[\rm (a)]
\item The following are equivalent:
\begin{enumerate}[\rm (i)]
	\item $r(X)$ and $s(X)$ have the same Frobenius norm, for all $X\in\dom r\cap \dom s$;
	\item $r=\alpha+u$ and $s=\xi\alpha+\zeta u$ for some $\alpha\in\kk$, $u\in U$, and $\xi,\zeta\in\kk$ of modulus one such that either $\xi=\zeta$ or $u\in [U,U]$.
\end{enumerate}
\item The following are equivalent:
\begin{enumerate}[\rm (i')]
	\item $r(X)$ and $s(X)$ have the same operator norm, for all $X\in\dom r\cap \dom s$;
	\item the singular values of $r(X)$ and $s(X)$ coincide, for all $X\in\dom r\cap \dom s$;
	\item $s=\zeta r$ for some $\zeta\in\kk$ of modulus one.
\end{enumerate}
\end{enumerate}
\end{corollary}

\begin{proof}
If one of $r,s$ is constant, so is the other, since a noncommutative rational function with bounded norm on $\dom_n r$ is constant on $(\kk^{n\times n})^\gen$ by Liouville's theorem. Thus, it suffices to assume that $r,s\in U\setminus\kk$.

(a): Suppose (ii) holds. If $\xi=\zeta$, then (i) is immediate. If $u\in[U,U]$, then $\tr u(X)=0$
on $\dom u$, and so
\begin{align*}
\|s(X)\|_F^2&=
\tr\Bigl(
(\xi\alpha I+\zeta u(X))^* (\xi\alpha I+\zeta u(X))
\Bigr) \\
&=\tr\Bigl(
|\alpha|^2 I+u(X)^* u(X)
\Bigr) \\
&=\tr\Bigl(
(\alpha I+u(X))^* (\alpha I+u(X))
\Bigr)
=\|r(X)\|_F^2.
\end{align*}

Now suppose (i) holds. Then $\tr((rr^*-ss^*)(X,X^*))=0$ for all $X\in\dom r\cap\dom s$. By Zariski denseness of \eqref{e:dense}, we have $\tr((rr^*-ss^*)(X,Y))=0$ for all $X,Y\in\dom r\cap\dom s$. For every matrix $H$ and sufficiently small $\varepsilon\in\kk$,
\begin{align*}
0&=\tr\Big((rr^*-ss^*)\big(X,(I-\varepsilon H)Y(I-\varepsilon H)^{-1}\big)\Big) \\
&=\tr\Big(
r(X)(I-\varepsilon H)r^*(Y)(I-\varepsilon H)^{-1}
-s(X)(I-\varepsilon H)s^*(Y)(I-\varepsilon H)^{-1}
\Big).
\end{align*}
Expanding the right-hand side as a power series in $\varepsilon$ and taking the linear term gives
$$\tr\Big(\big((rr^*-r^*r-ss^*+s^*s)(X,Y)\big)H\Big)=0$$
for all $X,Y,H$, and so $rr^*-r^*r=ss^*-s^*s$. By Lemma \ref{l:norm}(a) we have $s=\zeta r+\lambda$ for $\lambda,\zeta\in\kk$ with $|\zeta|=1$. If $\lambda=0$, then (ii) holds. If $\lambda\neq0$, then $\tr((rr^*-ss^*)(X,Y))=0$ gives
$$\bar{\lambda}\zeta \tr r(X)
+\lambda\bar{\zeta} \tr r^*(Y)+n|\lambda|^2=0$$
for all $X,Y\in \dom_n r$ and $n\in\N$. Since $\lambda\neq0$, this is only possible if there is $\alpha \in\kk^\times$ such that $\tr r(X)= n\alpha$ for all $X\in \dom_n r$ and $n\in\N$. Then, $r-\alpha\in [U,U]$ by Theorem \ref{t:add}. Furthermore, $\bar{\lambda}\zeta \alpha+\lambda\bar{\zeta} \bar{\alpha}+|\lambda|^2=0$, so $\xi=\frac{\zeta\alpha +\lambda}{\alpha}$ satisfies
$|\xi|=1$, and (ii) holds.

(b): The implications (iii')$\implies$(ii')$\implies$(i') are clear.
The proof of implication (i')$\implies$(ii') translates almost verbatim from the polynomial case established in \cite[Theorem 6.3]{SV}, so we only provide its sketch.
If the operator norms of $r$ and $s$ agree, then the largest eigenvalues of $rr^*$ and $ss^*$ agree on evaluations of the form $(X,X^*)$. Thus, the characteristic polynomials of $(rr^*)(X,X^*)$ and $(ss^*)(X,X^*)$ share a root, for every $X\in\dom r\cap\dom s$. By Zariski denseness of \eqref{e:dense} and the resultant test for common roots, it follows that characteristic polynomials of $(rr^*)(X,Y)$ and $(ss^*)(X,Y)$ share a root for all $X,Y\in\dom r\cap\dom s$. 
For sufficiently large prime size $n$, the characteristic polynomials of $rr^*$ and $ss^*$ evaluated at the $2\gen$-tuple of generic $n\times n$ matrices are irreducible \cite[proof of Lemma 2.6]{BV}. Since they share a root, they coincide. Thus, $(rr^*)(X,Y)$ and $(ss^*)(X,Y)$ have the same eigenvalues, so $r(X)$ and $s(X)$ have the same singular values for all $X\in \dom r\cap\dom s$.

Finally, assume (ii'). Then $\tr(((rr^*)^{\pm1}-(ss^*)^{\pm1})(X,X^*))=0$ for all $X$. As in the proof of (i)$\implies$(ii) in (a), we conclude that $(rr^*)^{\pm1}-(r^*r)^{\pm1}=(ss^*)^{\pm1}-(s^*s)^{\pm1}$. Then (iii') holds by Lemma \ref{l:norm}(b).
\end{proof}

\section{Multiplicative commutators}\label{s:mult}

We establish a multiplicative counterpart of
Theorem~\ref{t:add}, characterising products of multiplicative
commutators by the determinants of their matrix evaluations
(Theorem~\ref{t:mult}). As a consequence, we show that two
noncommutative rational functions have the same spectra at
every common matrix evaluation if and only if they are conjugate
in the free skew field (Corollary~\ref{c:isospec}).

Let $A$ be an $m\times (m+1)$ matrix over $R$ such that there is a unique $u\in U^{m+1}$ with the first entry $1$ that satisfies $Au=0$. If $r\in U$ is the last entry of $u$, we say that such an $A$ is \emph{admissible for $r$}; by \cite[Section 7.1]{Coh}, there exist admissible matrices for every element of $U$. If we single out the first and the last column of $A$,
$$A=(A_0\ A_\bullet\ A_\infty),\qquad (A_0\ A_\bullet),(A_\bullet\ A_\infty)\in R^{m\times m},$$
then $(A_0\ A_\bullet)$ and $(A_\bullet\ A_\infty)$ are viewed as the numerator and the denominator, respectively, of the representation $Au=0$ for $r$. Note that uniqueness of $u$ is equivalent to the denominator being invertible over $U$, or equivalently, that the denominator is a full matrix over $R$. 
If $r\in U^\times$, then $(A_\infty\ A_\bullet\ A_0)$ is admissible for $r^{-1}$, so $(A_0\ A_\bullet)$ is also full over $R$. By Cramer's rule \cite[Proposition 7.1.5]{Coh},
$$(A_\bullet\ -A_0)=(A_\bullet\ A_\infty)\begin{pmatrix}
I & u_\bullet \\ 0 & r
\end{pmatrix},$$
where $u_\bullet$ is the middle vector of $u$, 
i.e., $u^\intercal=(1\ u_\bullet^\intercal \ r)$.

\def\Full{\operatorname{Full}}

A full matrix $X$ over $R$ is \emph{irreducible} if it is not invertible over $R$ and is not a product of two non-invertible full matrices over $R$. 
Full matrices $X,Y$ over $R$ are \emph{stably associated} \cite[Section 0.5]{Coh} if $P(X\oplus I)=(Y\oplus I)Q$ for some invertible matrices $P,Q$ over $R$. 

Let $\Full_m(R)$ denote the set of all full $m\times m$ matrices over $R$. 
The embeddings $\Full_m(R)\hookrightarrow \Full_{m+1}(R)$ given by $X\mapsto X\oplus 1$ give rise to the direct limit $\Full(R)=\varinjlim \Full_m(R)$. 
Stable association and irreducibility are compatible with these embeddings, so stable association gives an equivalence relation on $\Full(R)$, and it preserves irreducible elements.

The divisor group $\operatorname{D}(R)$ \cite[Section 7.10]{Coh} is the free abelian group generated by stable association classes of irreducible matrices over $R$. Since factorisation into irreducible matrices is unique up to permutation and stable association of the irreducible factors \cite[Proposition 3.2.9]{Coh}, there is a monoid homomorphism
$$\phi:\Full(R)\longrightarrow \operatorname{D}(R),\qquad \phi(X_1\cdots X_\ell)=[X_1]+\cdots+[X_\ell],$$
where $[X]\in \operatorname{D}(R)$ denotes the class of irreducible $X\in\Full(R)$. 
By \cite[Theorem 7.10.2]{Coh}, this map gives rise to a well-defined group homomorphism 
\begin{align*}
\varphi:U^\times_{\operatorname{ab}}\longrightarrow \operatorname{D}(R),\qquad 
 \varphi([r])&=\phi(A_0\ A_\bullet)-\phi(-A_\infty\ A_\bullet) \\
& \text{if }(A_0\ A_\bullet\ A_\infty)\text{ is admissible for }r,
\end{align*}
where $[r]\in U^\times_{\operatorname{ab}}\coloneq U^\times/(U^\times)'$ denotes the class of $r\in U^\times$. The kernel of $\varphi$ is $\kk^\times$ by \cite[Theorem 7.10.9]{Coh}. 
We collect all these facts in the following statement.

\begin{proposition}\label{p:coh}
There is a group homomorphism $\varphi: U^\times_{\operatorname{ab}}\to \operatorname{D}(R)$ with kernel $\kk^\times$, such that for every $r\in U^\times$ there are $X,Y\in\Full(R)$ such that
$$
\varphi([r])=\phi(X)-\phi(Y),\qquad
\begin{pmatrix}
r&0\\ \star & I
\end{pmatrix}=Y^{-1}X.
$$
\end{proposition}

Next, we record that the commutator subgroup of $U^\times$ behaves well under scalar extensions.

\begin{lemma}\label{l:ac}
Let $\overline{\kk}$ be the algebraic closure of $\kk$. Then, $U^\times\cap (U_{\overline{\kk}}^\times)'=(U^\times)'$. 
\end{lemma}

\begin{proof}
Let $r\in U^\times\cap (U_{\overline{\kk}}^\times)'$. Since only finitely many elements of $\overline{\kk}$ appear in the presentation of $r$ as a product of multiplicative commutators, we have $r\in  U^\times\cap (U_K^\times)'$ for some finite extension $K$ of $\kk$. 
Fix an embedding $K\hookrightarrow \kk^{n\times n}$, and consider the group homomorphism
\begin{equation}\label{e:dieudonne}
U_K^\times\cong (K\otimes_\kk U)^\times\to \GL_n(U)\to U^\times_{\operatorname{ab}},
\end{equation}
where the last map is the Dieudonn\'e determinant. Then, \eqref{e:dieudonne} maps $r$ to $[r^n]$. Since the Dieudonn\'e determinant maps commutators to the identity, we have $[r]^n=[1]$. Since $U^\times_{\operatorname{ab}}/\kk^\times$ is a free abelian group by Proposition \ref{p:coh}, there is $\alpha\in\kk^\times$ such that $\alpha r\in (U^\times)'$. Then, $\alpha\in (U_K^\times)'\cap K^\times=\{1\}$, and therefore $r\in (U^\times)'$.
\end{proof}

Combining Proposition \ref{p:coh} with properties of determinants of noncommutative polynomials \cite{HKV} allows us to characterise noncommutative rational functions with determinant one, as follows.

\begin{theorem}\label{t:mult}
Let $r\in U^\times$. Then, $\det r= 1$ on $\dom r$ if and only if $r \in (U^\times)'$.
\end{theorem}

\begin{proof}
$(\impliedby)$ If $r=(a_1b_1a_1^{-1}b_1^{-1})\cdots (a_\ell b_\ell a_\ell^{-1}b_\ell^{-1})$ for $a_i,b_i\in U^\times$, 
multiplicativity of the determinant gives $\det r= 1$ on $\bigcap_i\dom a_i\cap \dom a_i^{-1}\cap \dom b_i\cap \dom b_i^{-1}$, and then $\det r= 1$ on $\dom r$ by eventual Zariski denseness.

$(\implies)$ Let $r\notin (U^\times)'$. 
By Lemma \ref{l:ac}, it suffices to assume that $\kk$ is algebraically closed. 
If $[r]=[\alpha]$ for $\alpha\in\kk^\times\setminus\{1\}$, then $\det r=\alpha^n$ on $\dom_n r$ for all $n$, so $\det r$ is not equal to $1$ on the whole $\dom r$. 
If $[r]\notin \kk^\times\subset U^\times_{\operatorname{ab}}$, let $X,Y$ be as in Proposition \ref{p:coh}. Then,
$$\varphi([r])=\phi(X)-\phi(Y)=\sum_{i=1}^N n_i[F_i]$$
for some $N\in\N$, nonzero $n_i\in\Z$, and irreducible matrices $F_i$ over $R$ from pairwise distinct stable association classes. Let $\Omega_n$ denote the $\gen$-tuple of generic $n\times n$ matrices (i.e., its $\gen n^2$ entries are independent commuting indeterminates). For all sufficiently large $n\in\N$, 
there is $\alpha_n\in\kk^\times$ such that
$$\det r(\Omega_n)=
\det (Y(\Omega_n))^{-1}
\det (X(\Omega_n))
=\alpha_n\prod_{i=1}^N \det (F_i(\Omega_n))^{n_i}$$
by the definition of stable association, because $\det P(\Omega_n)\in \kk^\times$ for every invertible matrix $P$ over $R$. For sufficiently large $n$, the commutative polynomials $\det (F_1(\Omega_n)),\dots,\det (F_N(\Omega_n))$ are irreducible by \cite[Theorem 2.9]{HKV}, and pairwise non-associated by \cite[Corollary 2.13]{HKV}. Thus, $\det r(\Omega_n)$ is not constant for some $n$.
\end{proof}

As a consequence, we show that noncommutative rational functions with matching spectra (in the algebraic closure of $\kk$) are conjugate (see \cite[Theorem 1.2]{SV} for a polynomial analogue).

\begin{corollary}\label{c:isospec}
Let $r,s\in U$. The spectra of $r(X)$ and $s(X)$ coincide for all $X\in\dom r\cap \dom s$ if and only if there is $u\in U^\times$ such that $s=uru^{-1}$.
\end{corollary}

\begin{proof}
Without loss of generality, $r,s\not\in\kk$.
Only the forward implication is nontrivial. Suppose that the spectra of $r(X)$ and $s(X)$ coincide on $\dom r\cap \dom s$. For a generic matrix tuple $X$ (of sufficiently large dimension), all the eigenvalues of $r(X)$ and $s(X)$ are simple by \cite[Theorem 2.8]{BV}. Since the characteristic polynomial of a matrix with only simple eigenvalues is determined by these eigenvalues, it follows by Zariski density that $\det (\lambda I-r(X))=\det (\lambda I-s(X))$ for all $\lambda\in\kk$ and $X\in\dom r\cap\dom s$. Consequently, $(t-r)(t-s)^{-1}\in U_{\kk(t)}$ has constant determinant one, and so $(t-r)(t-s)^{-1}\in (U_{\kk(t)}^\times)'$ by Theorem \ref{t:mult}.

The ring $U[t]$ is a unique factorisation domain \cite[Theorems 1.3.2 and 1.3.5]{Coh}: if $p_1\cdots p_\ell=q_1\cdots q_m$ are two factorisations into irreducible elements of $U[t]$, then $\ell=m$ and there is a permutation $\pi$ of $\{1,\dots,\ell\}$ such that every $q_i$ is \emph{similar} to $p_{\pi(i)}$, meaning $U[t]/q_iU[t]\cong U[t]/p_{\pi(i)}U[t]$ as right $U[t]$-modules. Thus, every irreducible $p\in U[t]$ determines a monoid homomorphism $\nu_p:(U[t]\setminus\{0\},\cdot)\to(\N_0,+)$, where $\nu_p(f)$ counts the number of irreducible factors of $f$ that are similar to $p$. Since $U(t)$ is the classical ring of quotients of $U[t]$, this map extends to a unique group homomorphism $\nu_p:U(t)^\times\to\Z$ satisfying $\nu_p(fg^{-1})=\nu_p(f)-\nu_p(g)$.

Note that $t-r,t-s\in U[t]$ are irreducible. Since $(t-r)(t-s)^{-1}\in (U(t)^\times)'$ and the image of $\nu_{t-r}$ is abelian, it follows that
$$\nu_{t-r}\bigl( (t-r)(t-s)^{-1} \bigr)=0,$$
so $t-r$ and $t-s$ are similar. A right $U[t]$-module isomorphism
$$U[t]/(t-r)U[t]\to U[t]/(t-s)U[t]$$
maps $1+(t-r)U[t]\mapsto u+(t-s)U[t]$ for some $u\in U^\times$ satisfying $u(t-r)\in (t-s)U[t]$. Comparing constant and linear coefficients shows that $ur=su$.
\end{proof}

\begin{remark}\label{r:semaev}
Theorems \ref{t:add} and \ref{t:mult} describe noncommutative rational functions whose pointwise eigenvalues add to zero or multiply to one. One can also easily find $r\in U$ whose pointwise eigenvalues satisfy other related equations; e.g., if $r$ is an inverse of an element in $[U,U]$, then for any invertible value of $r$, the reciprocals of its eigenvalues add to 0. What is harder to see is that there are eigenvalue relations that are not immediately derived from additive or multiplicative commutators.

For example, let $\kk=\C$ and let $p\in\kk[x]$ be a monic cubic polynomial with simple roots. Then $\mathscr C=\{(x,y)\in\kk^2: y^2=p(x)\}$ is (an affine chart of) an elliptic curve. By \cite[Corollary 5.3]{DV} there exist nonconstant commuting $r,s\in U$ such that $s^2=p(r)$. By the proof of \cite[Proposition 5.5]{DV}, the eigenvalues of $r(X)$ for $X\in\dom_n r$ ($n\ge2)$ are the $x$-coordinates of $n$ points on $\mathscr C$ whose sum, in the group law on $\mathscr C$, is the point at infinity. That is, the eigenvalues of $r(X)$ determine a zero of the Semaev polynomial $f_n$, which is irreducible and of total degree $(n-1)2^{n-2}$ \cite[Theorem 1]{Semaev}. 
This poses limitations on expressing the relation $f_n=0$ for eigenvalues of $r$ on $\dom_n r$ in terms of traces or determinants in an $n$-independent way. More precisely, there are no univariate rational functions $\alpha_1,\dots,\alpha_s,\beta_1,\dots,\beta_t$ and $(s+t)$-variate rational function $\rho$ such that $f_n$ evaluated on eigenvalues of a generic $n\times n$ matrix $Z$ agrees with 
$\rho(\tr\alpha_1(Z),\dots,\tr\alpha_s(Z),
\det\beta_1(Z),\dots,\det\beta_t(Z))$, for all $n$.
\end{remark}

\section{Products of additive commutators}\label{s:one}

Gardella and Thiel asked whether every element of a
noncommutative division ring is a product of two additive
commutators \cite[Question~5.9]{GT}. This has been
established for finite-dimensional noncommutative division algebras by
Bre\v{s}ar, Jang and Robert \cite{BJR}.
In contrast, we show that the identity in the free skew
field cannot be expressed as a nonempty product of additive
commutators (Corollary~\ref{c:one}), answering
\cite[Question~5.9]{GT} in the negative. Our proof uses a
valuation on an extension of the free skew field and a
corresponding obstruction in a free Poisson field.
Note that every element of $U^{n\times n}$ for $n\ge 2$ is a
product of two additive commutators by
\cite[Theorem 4.4]{BGT}.

If $\gen>1$, then every element in $U$ is a sum of two products of two additive commutators by \cite[Proposition 5.8]{GT}. For instance, $1$ equals
$$
\Bigl[\bigl[[x_1^2,x_2],[x_1,x_2]\bigr]^{-1}[x_1^2,x_2][x_1,x_2]^{-1},[x_1,x_2]\Bigr]
[x_1,x_2]+
\Bigl[[x_1,x_2],\bigl[[x_1^2,x_2],[x_1,x_2]\bigr]^{-1}\Bigr]
[x_1^2,x_2].
$$
However, by Corollary \ref{c:one} below, 1 and the generators $x_i$ are not nonempty products of additive commutators.  
For completeness, let us mention that the $x_i$ are not sums of multiplicative commutators \cite[Exercise 7.10.8]{Coh}.

Let $L$ be the free Lie $\kk$-algebra generated by $x_1,\dots,x_\gen$ \cite[Section 0.2]{Reu}. Then, $L$ is graded by Lie word length; let $\deg_L$ be the resulting degree function on $L$, and let $(e_i)_i$ be a homogeneous basis of $L$. 
Its corresponding free Poisson field \cite[Section 2]{MS} is 
the rational field $P=\kk(z_i\colon i)$ with the Poisson bracket $\{\cdot,\cdot\}$ satisfying
$$\{z_i,z_j\}=\sum_k \gamma_{ijk}z_k\qquad\text{when}\qquad [e_i,e_j]=\sum_k \gamma_{ijk}e_k.$$
For $n\in\N$ denote $P_n=\kk(z_i\colon \deg_L e_i\le n)$.

\begin{lemma}\label{l:poisson}
No element of $P_1^\times$ is a nonempty product of Poisson brackets
in $P$.
\end{lemma}

\begin{proof}
First, we prove an auxiliary assertion for any $f,g\in P$ with $\{f,g\}\neq0$. Namely, if $n\in\N$ is minimal such that $f,g\in P_n$, then
\begin{equation}\label{e:deg}
\{f,g\}\in P_n[z_i\colon \deg_L e_i>n]\setminus P_n.
\end{equation}
Indeed, consider $F=\sum_i\frac{\partial f}{\partial z_i}e_i$ and $G=\sum_i\frac{\partial g}{\partial z_i}e_i$ in $P_n\otimes_\kk L$. The $P_n$-linear map $P_n\otimes_\kk L\to P$ determined by $e_i\mapsto z_i$ maps $[F,G]$ to
$$\{f,g\}=\sum_{i,j} \frac{\partial f}{\partial z_i}\frac{\partial g}{\partial z_j}\{z_i,z_j\}
\in P_n[z_i\colon \deg_L e_i>n].$$
Next, $\{f,g\}\neq0$ implies $[F,G]\neq0$, and we claim that $\deg_L[F,G]>n$.
Indeed, let $F_p,G_q$ be the leading homogeneous
components of $F,G$ (so $\max\{p,q\}=n$ by the assumption on $n$). If $[F_p,G_q]\neq0$, then
$\deg_L[F,G]=p+q>n$.
Otherwise, $F_p,G_q$ are proportional by the
Shirshov--Witt theorem \cite[Theorem 2.5]{Reu},
so $p=q$. Replacing $G$ by $G-cF$ for a suitable $c\in P_n$ does not change $[F,G]$ and
lowers the degree of $G$. The leading components
now have distinct degrees, hence cannot commute
by the Shirshov--Witt theorem. Therefore
$\deg_L[F,G]>p=n$. 
Consequently, $\{f,g\}\notin P_n$.

Now suppose $\prod_{j=1}^\ell\{f_j,g_j\}\in P_1^\times$, and let $n\in\N$ be minimal such that
$f_j,g_j\in P_n$ for all $j$. Then, $\{f_j,g_j\} \in P_n[z_i\colon \deg_L e_i>n]$ for all $j$ 
and $\{f_j,g_j\}\notin P_n$ for some $j$ by \eqref{e:deg}.
Consequently, $\prod_{j=1}^\ell\{f_j,g_j\}\notin P_n$, a contradiction.
\end{proof}

\def\lc{\operatorname{lc}}

The universal enveloping algebra of $L$ is $R$. By \cite[Theorem 1 and Proposition 2.3]{Lic}, the free skew field $U$ embeds into a skew field $D$ with a discrete valuation $\nu:D\to\Z\cup\{\infty\}$ and a central element $t$ such that
$$\nu(\kk^\times)=\{0\},\qquad \nu(t)=1,\qquad \nu(e_i)=-1,\qquad
 \cO/t\cO\cong P,\qquad \overline{te_i}=z_i,
$$
where $\cO=\{u\in D:\nu(u)\geq0\}$ (here, $\overline{a}\in P$ denotes the coset of $a\in\cO$). 
Observe that
\begin{equation}\label{eq:comm-bound}
 \nu([a,b])\geq\nu(a)+\nu(b)+1
 \qquad\text{for all }a,b\in D^\times.
\end{equation}
Indeed, since the residue field is commutative, the commutator of $t^{-\nu(a)}a,t^{-\nu(b)}b\in\cO$ is not invertible in $\cO$, so its value is strictly positive, implying \eqref{eq:comm-bound}.
Since
\[
 \overline{t^{-1}[te_i,te_j]}
 =\overline{t[e_i,e_j]}
 =\sum_k \gamma_{ijk}%
 z_k,
\]
we have $\overline{t^{-1}[a,b]}=\{\overline{a},\overline{b}\}$ for $a,b\in\cO$ 
in terms of the Poisson bracket on $P$.
For $a\in D^\times$, define
\[
 \lc(a)\coloneq\overline{t^{-\nu(a)}a}\in P^\times.
\]
Then,
\begin{equation}\label{eq:lc}
\lc(ab)=\lc(a)\lc(b),\qquad
 \lc(t^n a)=\lc(a)\quad\text{for }n\in\Z.
\end{equation}

\begin{lemma}\label{l:lead}
If $a,b\in D$ and $[a,b]\ne0$, then $\lc([a,b])\in P$ is a Poisson bracket.
\end{lemma}

\begin{proof}
After scaling with suitable powers of $t$, it suffices to assume $\nu(a)=\nu(b)=0$ by \eqref{eq:lc}. 
We prove the statement by induction on $m=\nu([a,b])$. 
If $m=1$, then $\lc([a,b])=\overline{t^{-1}[a,b]}=\{\overline{a},\overline{b}\}$ is a Poisson bracket. 

Now suppose $m>1$. Then, $t^{-1}[a,b]$ is not invertible in $\cO$, 
so $\{\overline{a},\overline{b}\}=\overline{t^{-1}[a,b]}=0$. By \cite[Theorem 1]{MS}, $\overline{a},\overline{b}\in P$ are algebraically dependent, so there is $0\neq p\in \kk[x,y]$ of minimal degree such that $p(\overline{a},\overline{b})=0$. Without loss of generality, $\frac{\partial p}{\partial y}(\overline{a},\overline{b})\neq0$. Write $p=\sum_jp_j(x)y^j$, and set
$$A=\sum_j p_j(a)b^j,\qquad H=\sum_j jp_j(a)b^{j-1}.$$
Then $A\in t\cO$ and $H\in\cO^\times$ by the choice of $p$. 
For every $k\in\N$,
\[
 [a,b^k]=\sum_{j=0}^{k-1}b^j [a,b] b^{k-1-j}
 \equiv [a,b]\cdot k b^{k-1}\pmod{t^{m+1}\cO}
\]
because $[a,b]$ commutes with elements of $\cO$ modulo $t^{m+1}\cO$ by \eqref{eq:comm-bound}. Hence,
\[
 [a,A]\equiv [a,b]H\pmod{t^{m+1}\cO},
 \qquad \nu([a,H])\geq m
\]
and consequently,
\[
 [a,AH^{-1}]=[a,A]H^{-1}-AH^{-1}[a,H]H^{-1} \equiv [a,b]\pmod{t^{m+1}\cO}.
\]
In particular, $\nu([a,AH^{-1}])=m$ and
$\lc([a,AH^{-1}])=\lc([a,b])$. Put $q=\nu(AH^{-1})\geq1$ and $b'=t^{-q}AH^{-1}$.
Since $\nu(b')=0$ and
\[
 1\le \nu([a,b'])=m-q<m
\]
by \eqref{eq:comm-bound}, $\lc([a,b])=\lc([a,b'])$ is a Poisson bracket by the induction hypothesis.
\end{proof}

\begin{corollary}\label{c:one}
A nonzero element of $\operatorname{span}_\kk\{1,x_1,\dots,x_\gen\}$ is not a nonempty product of additive commutators in $U$.
\end{corollary}

\begin{proof}
Let $a$ be a nonzero linear combination of $1$ and the generators $x_i$. As an element the universal enveloping algebra of $L$, $a$ is a linear combination of 1 and the $e_i$ with $\deg_L e_i=1$, so $\lc(a)\in P_1^\times$.
If $a=[b_1,c_1]\cdots [b_\ell,c_\ell]$, then $\lc(a)=\prod_{j=1}^\ell \lc([b_j,c_j])$ by \eqref{eq:lc}. Hence, $\lc(a)\in P_1^\times$ is a product of Poisson brackets by Lemma \ref{l:lead}, which contradicts Lemma \ref{l:poisson}.
\end{proof}

\end{document}